\documentclass[a4paper,11pt]{article}
\usepackage[english]{babel}							
\usepackage{amssymb,amsmath,cancel,mathrsfs} 

\usepackage{epsfig}
\usepackage{slashed}

\def\beq{\arraycolsep1pt\begin{eqnarray*}}
\def\eeq{\end{eqnarray*}}

\newcommand\R[0]{\mathbb{R}}
\newcommand\N[0]{\mathbb{N}}

\newcommand\ds[0]{\displaystyle}
\newcommand\ee{\varepsilon}
\newcommand\pscal[2]{\left\langle \, #1 \, , \, #2 \, \right\rangle}
\newcommand\aand[0]{\text{ and }}
\newcommand\forevery[0]{\text{ for every }}
\newcommand\cvd[0]{\hfill$\blacksquare$}
\newcommand\proof[0]{{\it Proof. }}

\begingroup
\newtheorem{theorem}{Theorem}[section]

\newtheorem{lemma}[theorem]{Lemma}

\newtheorem{remark}[theorem]{Remark}
\newtheorem{proposition}[theorem]{Proposition}
\newtheorem{example}[theorem]{Example}
\newtheorem{claim}[theorem]{Claim}
\endgroup

\title{Double resonance for one-sided superlinear or singular nonlinearities} 

\author{Andrea Sfecci}

\date{}

\begin{document}

\maketitle

\begin{abstract}
We deal with the problem of existence of periodic solutions for the scalar differential equation $x''+ f(t,x)=0$ when the {\em asymmetric} nonlinearity satisfies a {\em one-sided} superlinear growth at infinity. The nonlinearity is asked to be {\em next to} resonance and a Landesman-Lazer type of condition will be introduced in order to obtain a positive answer. Moreover we provide also the corresponding result for equations with a singularity and asymptotically linear growth at infinity, showing a further application to radially symmetric systems.

\end{abstract}

\section{Introduction}\label{intro}

In this paper we are going to study different types of scalar second order differential equations. We are interested in nonlinearities which, roughly speaking, are {\em next to resonance}. We will provide different results of existence of periodic solutions extending some previous theorems well-known in literature treating the case of {\em nonresonant} nonlinearities. We will first focus our attention on nonlinearities defined on the whole real line, in particular, we will start looking for periodic solutions of the scalar differential equation
\begin{equation}\label{nleq}
x''+ f(t,x)=0 \,,
\end{equation}
where $f: \R \times \R \to \R$ is a continuous function which is $T$-periodic in the first variable. Then, in Section~\ref{secrad}, we will show how such a result can be adapted to the case of nonlinearities with a singularity. Finally, we will provide a further application to radially symmetric systems.

\medbreak

It is well-known by classical results~\cite{Dancer,Fucik,MawhinSurvey} that the asymmetric oscillator
$$
x''+\mu x^+ -\nu x^- =0
$$
has nontrivial solutions if the couple $(\mu,\nu)$ belongs to the so-called Dancer-Fucik spectrum
$$
\Sigma= \bigcup_{j\in \N} \mathscr C_j\,,
$$
where
$$
\mathscr C_0 = \left\{ (\mu,\nu) \in \R^2 \,:\,  \mu\geq0, \nu\geq0 \text{ such that } \mu\,\nu=0\right\}
$$
and
$$
\mathscr C_j = \left\{ (\mu,\nu) \in \R^2 \,:\, \mu\geq0, \nu\geq0 \text{ such that } \frac{\pi}{T} \left(\frac{1}{\sqrt{\mu}} + \frac{1}{\sqrt{\nu}}\right)= \frac 1j \right\}.
$$
In particular, it consists of the two positive semi-axes $\mathscr C_0$ and of an infinite number of curves $\mathscr C_j$ having a vertical asymptote $\mu=\mu_j$ and an horizontal one $\nu=\nu_j$ with $\mu_j=\nu_j=(j\pi/T)^2$.

The study of {\em asymmetric nonlinearities} $f$ satisfying
$$
\nu_\downarrow \leq \liminf_{x\to-\infty} \frac{f(t,x)}{x} \leq  \limsup_{x\to-\infty} \frac{f(t,x)}{x}\leq \nu_\uparrow\,,
$$
$$
\mu_\downarrow \leq \liminf_{x\to+\infty} \frac{f(t,x)}{x} \leq  \limsup_{x\to+\infty} \frac{f(t,x)}{x}\leq \mu_\uparrow\,,
$$
for some suitable constants in $[0,+\infty]$,
providing the existence of periodic solutions to equation~\eqref{nleq} presents a wide literature (see e.g.~\cite{CMZ,Dancer,DI,Fabry,FH,FS1,Fucik} or the survey~\cite{MawhinSurvey} and the references therein).
The existence is strictly related to the position of the rectangle $W =[\mu_\downarrow,\mu_\uparrow]\times[\nu_\downarrow,\nu_\uparrow]$ with respect to the set $\Sigma$:
\begin{itemize}
\item if $W\cap \Sigma = \varnothing$ ({\em non-resonance}) we have the existence of at least one periodic solution  (cf.~\cite{CMZ,DI,FH}),
\item if $W$ is bounded and $W\cap \Sigma =\{(\mu_\downarrow,\nu_\downarrow)\}$ or $W\cap \Sigma =\{(\mu_\uparrow,\nu_\uparrow)\}$ ({\em simple resonance}) the existence of a periodic solution can be obtained by the introduction of a Landesman-Lazer type of condition (cf.~\cite{Fabry}) or of an Ahmad-Lazer-Paul type of condition (cf.~\cite{J-ALP}),
\item if $W$ is bounded and $W\cap \Sigma =\{(\mu_\downarrow,\nu_\downarrow),(\mu_\uparrow,\nu_\uparrow)\}$ ({\em double resonance}) the existence of a periodic solution can be obtained by the introduction of a {\em double} Landesman-Lazer type of condition (cf.~\cite{Fabry,FondaPHH,FG,FM})\,.
\end{itemize}
In this paper we will present an existence result of periodic solutions for the {\em double resonance} case in which $W$ is unbounded. Such a situation presents three possible interpretations: nonlinearities with one-sided superlinear growth, nonlinearities with a singularity and scalar equations with impacts. In this paper we will present the first two situations, the last one has been considered by the author recently in \cite{S10}. 
We are going to treat nonlinearities satisfying the following asymptotic asymmetric behavior.

\begin{itemize}
\item[{\bf (A)}] {\sl Assume
\begin{equation}\label{as-}
\lim_{x\to - \infty} \frac{f(t,x)}{x} = + \infty\,
\end{equation}
and that there exists a constant $c>0$ and an integer $N>0$ such that
\begin{equation}\label{as+}
\mu_N x - c \leq f(t,x) \leq \mu_{N+1} x + c\,,
\end{equation}
for every $x>0$ and every $t\in[0,T]$, where $\mu_j=(j\pi/T)^2$.
}
\end{itemize}

\noindent Notice that the specular case can also be considered as well.
The case of a nonlinearity satisfying a  {\em nonresonant} one-sided superlinear growth was studied e.g. in~\cite{FH}, but to the best of our knowledge an existence result for nonlinearities satisfying a {\em double} resonance condition has not been provided yet.

\medbreak

We are now ready to state the first of the main results of this paper. We address the reader to Sections~\ref{secrad} for corresponding theorems related to scalar equations with a singularity and to Section~\ref{final} for some applications to radially symmetric systems.

\begin{theorem}\label{main}
Assume {\bf(A)} and the Landesman-Lazer conditions
\begin{equation}\label{LLcond1}
\int_0^T \liminf_{x\to+\infty} (f(t,x)-\mu_{N} x ) \phi_N(t+\tau) \, dt > 0\,,
\end{equation}
\begin{equation}\label{LLcond2}
\int_0^T \limsup_{x\to+\infty} (f(t,x)-\mu_{N+1} x ) \phi_{N+1}(t+\tau) \, dt < 0\,,
\end{equation}
where $\phi_j$ is defined as
$$
\phi_j (t)=
\begin{cases}
\sin (\sqrt{\mu_j} t) & t\in[0, T/j]\\
0 & t\in[T/j, T] \,.
\end{cases}
$$
Then, equation $x''+f(t,x)=0$ has at least one $T$-periodic solution.
\end{theorem}

It is possible to relax the Landesman-Lazer conditions~\eqref{LLcond1} and~\eqref{LLcond2} in the previous theorem, requiring that nonlinearity $f$ satisfies the following additional hypothesis. 

%
%
%
%
\begin{itemize}
\item[{\bf (H)}] {\sl For every $\tau\in[0,T]$ and for every $\zeta>0$, consider the set $\mathcal I(\tau,\zeta)=[\tau-\zeta,\tau+\zeta]$ and the functions
$$
f_{1,\zeta,\tau}(x) = \min_{t\in\mathcal I(\tau,\zeta)} f(t,x) 
\qquad
f_{2,\zeta,\tau}(x) = \max_{t\in\mathcal I(\tau,\zeta)} f(t,x) 
$$
with their primitives $F_{i,\zeta,\tau}(x)=\int_0^x f_{i,\zeta,\tau}(\xi) \,d\xi$. 
We assume that
$$
\lim_{\zeta \to 0} \left( \lim_{x\to-\infty} \frac{F_{2,\zeta,\tau}(x)}{F_{1,\zeta,\tau}(x)} \right) = 1
$$
uniformly in $\tau\in[0,T]$.
}
\end{itemize}

The variant of Theorem~\ref{main} is thus given by the next one.

\begin{theorem}\label{main2}
Assume {\bf(A)}, {\bf(H)}
and the Landesman-Lazer conditions~\eqref{LLcond1} and~\eqref{LLcond2}
where $\phi_j$ is defined as
$$
\phi_j (t)= \left| \sin (\sqrt{\mu_j} t) \right| \,.
$$
Then, equation $x''+f(t,x)=0$ has at least one $T$-periodic solution.
\end{theorem}

\medbreak

Let us briefly explain the main differences between the two types of Landesman-Lazer conditions adopted in the previous theorems. The one involved in Theorem~\ref{main} is stronger than the one introduced in Theorem~\ref{main2}. In fact, it is easy to verify that the first implies the second. Hence, we can replace the stronger Landesman-Lazer conditions of Theorem~\ref{main}, with the weaker ones by introducing the additional assumption {\bf(H)}. Roughly speaking, it requires that the superlinear behavior of the nonlinearity at $-\infty$ is an {\em infinity of the same order} when $t$ varies as explained in the following example, where we show some nonlinearities satisfying (or not) such a hypothesis.

\begin{example}\label{ex1}
Suppose that there exists a function $h:\R\to \R$ satisfying 
$$
\lim\limits_{x\to-\infty} \frac{h(x)}{x} = +\infty\,,
$$
such that
$$
0< \liminf_{x\to-\infty} \frac {f(t,x)}{h(x)} \leq \limsup_{x\to-\infty} \frac {f(t,x)}{h(x)} < +\infty \,.
$$
Then, {\bf (H)} holds. As a particular situation, we can consider a nonlinearity $f$ which can be split (when $x<0$) as 
$f(t,x)=q(t)h(x)+p(t,x)$ with $q(t)>0$ and $\lim\limits_{x\to-\infty} \frac{p(t,x)}{h(x)}=0$
uniformely in $t$. 

As a direct example, {\bf (H)} holds for nonlinearities $f$ not depending on $t$ when $x<0$, or nonlinearities as $f(t,x) = (1+\sin^2(t)) x^5 + x^3$, or $f(t,x)= x^3 + \sin^2(t) x^2$.

Otherwise, for example, if $f(t,x) = x^3 + \sin^2(t) x^5$ when $x<0$, then $f$ does not satisfy {\bf (H)}.
\end{example}


\section{Proof of Theorems~\ref{main} and~\ref{main2}}\label{secproof}

By degree theoretic arguments, the proof consists in finding a common {\em a~priori bound} for all the $T$-periodic solutions of the differential equations
\begin{equation}\label{homotopy}
x'' + g_\lambda(t,x) = 0\,,
\end{equation}
where $\lambda\in[0,1]$ and
$$
g_\lambda(t,x) = \lambda f(t,x) + (1-\lambda)h(t,x) \,,
$$
with
$$
h(t,x)=
\begin{cases}
f(t,x) & x<-1\\
\mu x + x[\,\mu x - f(t,x)] & -1\leq x\leq 0\\
\mu x & x>0
\end{cases}
$$
defining $\mu=(\mu_N+\mu_{N+1})/2$.

In particular we will find a $R_{good}>0$ such that every $T$-periodic solution of~\eqref{homotopy} satisfies $x(t)^2+ x'(t)^2 < R_{good}^2$ for every $t$. In~1993, Fabry and Habets proved in~\cite{FH} that there exists at least one $T$-periodic solution to~\eqref{homotopy} with $\lambda=0$, by the use of degree arguments. In particular, they found a similar a priori bound $R_{FH}$ for all the solutions of
$$
x''+ \tilde \lambda h(t,x) + (1-\tilde\lambda) (\mu x^+ - \nu x^-) =0 \,,
$$
with $\tilde\lambda\in[0,1]$ and $(\mu,\nu)\notin\Sigma$, thus solving the case of nonresonant nonlinearities. Hence, simply asking $R_{good}>R_{FH}$, using Leray-Schauder degree theory, the proof of Theorems~\ref{main} and~\ref{main2} easily follows. In Section~\ref{prel}, we will provide some preliminary lemmas which make use of some phase-plane techniques, then in Section~\ref{priori} we will prove the existence of the common a priori bound.

\subsection{Some preliminary lemmas}\label{prel}

In this section we present some estimates on the behavior of the solutions to~\eqref{homotopy} provided by the use of some phase-plane techniques. We will not present all the proofs, we leave some of them to the reader as an exercise of mere computation, referring to other papers for comparisons. By the way, some of the statements are well-known in literature.

\medbreak

Let us set 
$$
f_1(x)=\min_{t\in[0,T]} f(t,x) \quad \text{and}\quad
f_2(x)=\max_{t\in[0,T]} f(t,x) \,,
$$
then, by~\eqref{as-}, there exists $d<0$ such that
$$
f_1(x) < f(t,x) < f_2(x) < 0
$$
for every $x<d$,
with $\lim_{x\to-\infty} f_2(x)/x=+\infty$.
Define the primitives $ F_i (x)= \int_d^x f_i(\xi) d\xi$. Notice that $F_1>F_2$ are decreasing functions when $x<d$.

\medbreak
For every solution $x$ of equation~\eqref{homotopy} the couple $(x,y)=(x,x')$ is a solution of the planar system
\begin{equation}\label{planar}
\begin{cases}
x'=y\\
-y'= g_\lambda(t,x) \,.
\end{cases}
\end{equation}

\noindent We will say that
\begin{equation}\label{Rlarge}
(x,y) \text{ is $R$-large, if } x^2(t) + y^2(t) >R^2 \text{ for every } t\in[0,T]\,,
\end{equation}
where $(x,y)$ is a solution of~\eqref{planar}.

We will also consider the parametrization of such solutions in polar coordinates
$$
\begin{cases}
x(t) = \rho(t) \cos \theta(t)\\
y(t) = \rho(t) \sin \theta(t)
\end{cases}
$$
where the angular velocity and the radial velocity of $0$-large solutions are given by
$$
- \theta'(t) = \frac{y^2(t)+x(t)g_\lambda(t,x(t))}{x^2(t)+y^2(t)}\,,
$$
$$
\rho'(t) = \frac{y(t)[x(t)-g_\lambda(t,x(t))]}{\sqrt{x^2(t)+y^2(t)}}\,.
$$

An easy computation gives us the following.
\begin{remark}\label{rotate+}
There exists $R_0>0$ such that every $R_0$-large solution of~\eqref{planar} rotates clockwise (i.e. $- \theta'(t)>0$).
\end{remark}

In what follows the constant $R_0$ will be enlarged in order to obtain some additional properties on $R_0$-large solutions.

\begin{figure}[h]
\centerline{\epsfig{file=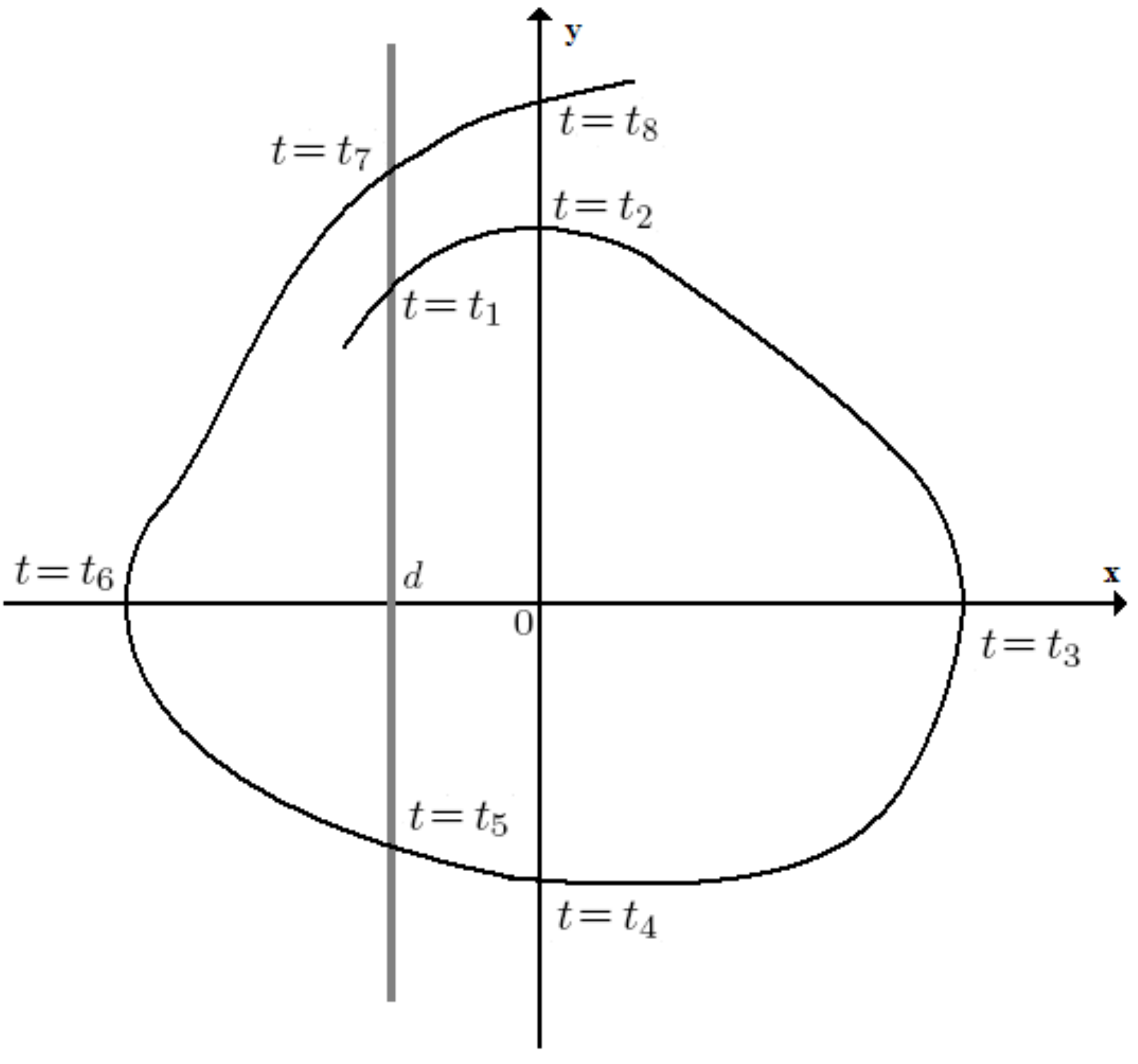, width = 9 cm}}
\caption{A $R_0$-large solution and the instants $t_1,\ldots,t_8$.}
\label{figlapt}
\end{figure}

Consider a $R_0$-large solution $(x,y)$, then there exist some instants $t_i$ (see Figure~\ref{figlapt}) such that
$$
\begin{array}{l}
x(t_1)=d, \quad y(t_1)>0\,,\\
x(t_2)=0, \quad y(t_2)>0\,,\\
x(t_3)>0, \quad y(t_3)=0\,,\\
x(t_4)=0, \quad y(t_4)<0\,,\\
x(t_5)=d, \quad y(t_5)<0\,,\\
x(t_6)<0, \quad y(t_6)=0\,,\\
x(t_7)=d, \quad y(t_7)>0 \,,\\
x(t_8)=0, \quad y(t_8)>0 \,.\\
\end{array}
$$

The following lemma can be obtained easily (see e.g.~\cite{FM,FS1} for details).

\begin{lemma}\label{nonreslemma}
For every $\ee>0$, it is possible to find $R_\ee>R_0$, such that every $R_\ee$-large solution of~\eqref{planar} satisfies
$$
t_5 - t_1\in\left( \frac{\pi}{\sqrt{\mu_{N+1}}} - \ee \,,\,  \frac{\pi}{\sqrt{\mu_{N}}} + \ee  \right)
=\left( \frac{T}{N+1} - \ee \,,\,   \frac{T}{N} + \ee  \right)
$$
and 
$t_7-t_5<\epsilon\,.$
\end{lemma}
%

Hence, we obtain easily the next one, choosing $\ee$ sufficiently small.

\begin{lemma}\label{Ngiri}
It is possible to find $R_0$, such that every $R_0$-large $T$-periodic solution of~\eqref{planar} performs exactly $N$ or $N+1$ rotations around the origin in the phase-plane.
\end{lemma}

With a similar reasoning, we can prove the following.
\begin{lemma}\label{halflap}
For every $\ee>0$, it is possible to find $R_\ee>R_0$, such that every $R_\ee$-large solution of~\eqref{planar} satisfies
$$
 \frac{T}{N+1} - \ee < t_4 - t_2 <   \frac{T}{N} + \ee \text{ and } t_8-t_4 < \ee \,.
$$
\end{lemma}

\medbreak

The following lemma gives us some informations on the dynamics when $x>d$.

\begin{lemma}\label{lapontheright}
It is possible to find $R_0$ large enough to have the existence of some positive constants $\theta_0$ and $\ell_0$ such that every $R_0$-large solution to~\eqref{planar}, when written in polar coordinates, satisfies
$$
-\vartheta'(t) > \omega_0 \quad\text{and}\quad |\rho'(t)| \leq \ell_0 \rho(t)
$$
when $x(t)>d$. So that, for $\kappa=\ell_0/\omega_0$,
$$
\left|\frac{d\rho}{d(-\theta)}\right| < \kappa \rho \,.
$$
\end{lemma}

We leave the proof to the reader as an exercise. We refer to~\cite{FS1,FS7} for details.

\begin{remark}\label{r1}
A direct consequence of the previous lemma is that
$$
e^{-\kappa\pi/2} \, x(t_3) \leq |y(t_j)| \leq e^{\kappa\pi/2} \, x(t_3)\,, \text{ with } j=2,4\,,
$$
if $(x,y)$ is $R_0$-large\,.
\end{remark}

\medbreak

Let us now focus our attention on the dynamics when $x<d$.
We are going to prove that there exists a functions $\mathcal T$ such that $y(t_7)< \mathcal T(y(t_5))$, thus permitting to find a second function $\mathcal L$ such that $y(t_8) \leq \mathcal L(y(t_2))$. We will have consequently a control on the behavior of solutions escaping from the origin.
 We start defining the energy functions
$$
H_i(x,y) = \frac 12 y^2 + F_i (x)\,, \quad i=1,2\,.
$$
Then, we have
\begin{equation}\label{en<1}
\frac {d}{dt} H_1(x(t),y(t))<0 \,\,\text{if } y(t)>0\,, \quad
\frac {d}{dt} H_1(x(t),y(t))>0 \,\,\text{if } y(t)<0\
\end{equation}
and
\begin{equation}\label{en<2}
\frac {d}{dt} H_2(x(t),y(t))<0 \,\,\text{if } y(t)<0\,, \quad
\frac {d}{dt} H_2(x(t),y(t))>0 \,\,\text{if } y(t)>0\,.
\end{equation}

\begin{figure}[h]
\centerline{\epsfig{file=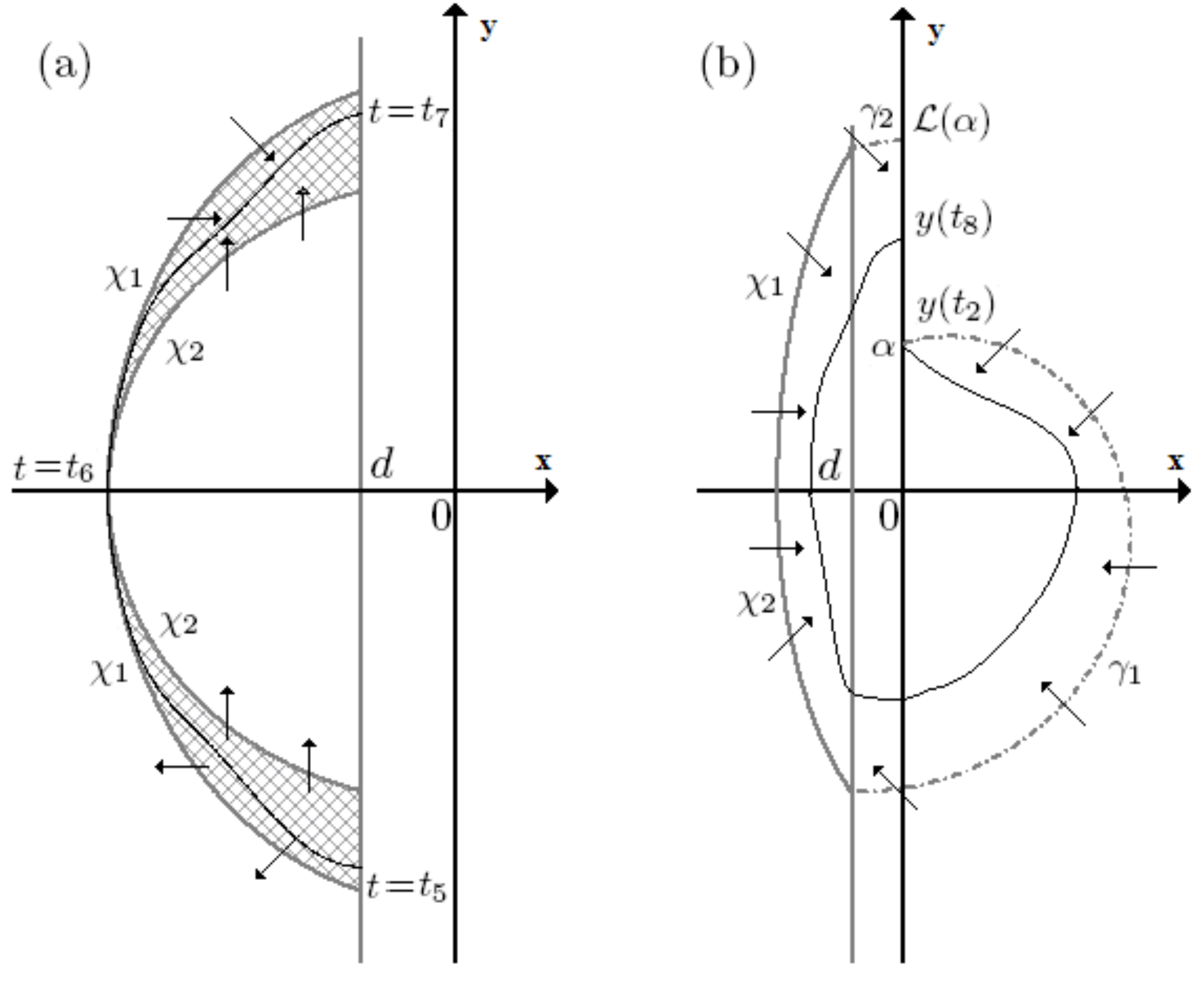, width = 12 cm}}
\caption{(a) The level subsets $\chi_1$ and $\chi_2$ of the energy functions $H_1$ and $H_2$ control the behavior of the solution $(x,y)$ of system~\eqref{planar} when $x<d$. (b) The solution, while it rotates around the origin, is {\em guided} by some curves, thus permitting to obtain the estimate in~\eqref{elle}. The curves $\gamma_1$ and $\gamma_2$ can be found by using the estimate in Lemma~\ref{lapontheright}.}
\label{figen}
\end{figure}

These functions give us a control on the behavior of the solutions, see Figure~\ref{figen}a. In particular, we have that
$$
-\sqrt{2F_1(x(t_6))} < y(t_5) < -\sqrt{2F_2(x(t_6))}\,,
$$
$$
\sqrt{2F_2(x(t_6))} < y(t_7) < \sqrt{2F_1(x(t_6))}\,,
$$
thus giving us $y(t_7)< \mathcal T(y(t_5))$, where
$$
\mathcal T(\upsilon) = \sqrt{2 F_1\left( F_2^{-1} (\upsilon^2/2) \right)} \geq \upsilon\,. 
$$
By the estimate in Lemma~\ref{lapontheright}, we have
$\sqrt{y(t_5)^2+d^2}<e^{\kappa\pi + a} \,y(t_2)$ and $y(t_8)< e^a \sqrt{y(t_7)^2+d^2}$,
where $a=\kappa\arcsin (d/R_0)$. Summing up, we obtain 
\begin{equation}\label{elle}
y(t_8) \leq \mathcal L(y(t_2)) \text{ with } \mathcal L(\upsilon)= e^a \sqrt{ \mathcal T^2 \left( \sqrt{ e^{2(\kappa\pi + a)} \, \upsilon^2 -d^2} \right) + d^2 }\,.
\end{equation}
 The same argument can be obtained by glueing together some {\em guiding curves} in the plane $(x,y)$ following an idea introduced in~\cite{FS1} and developed in~\cite{FS3,FS7,S5} in different situations. Figure~\ref{figen}b illustrates this idea.

\medbreak

\begin{figure}[h]
\centerline{\epsfig{file=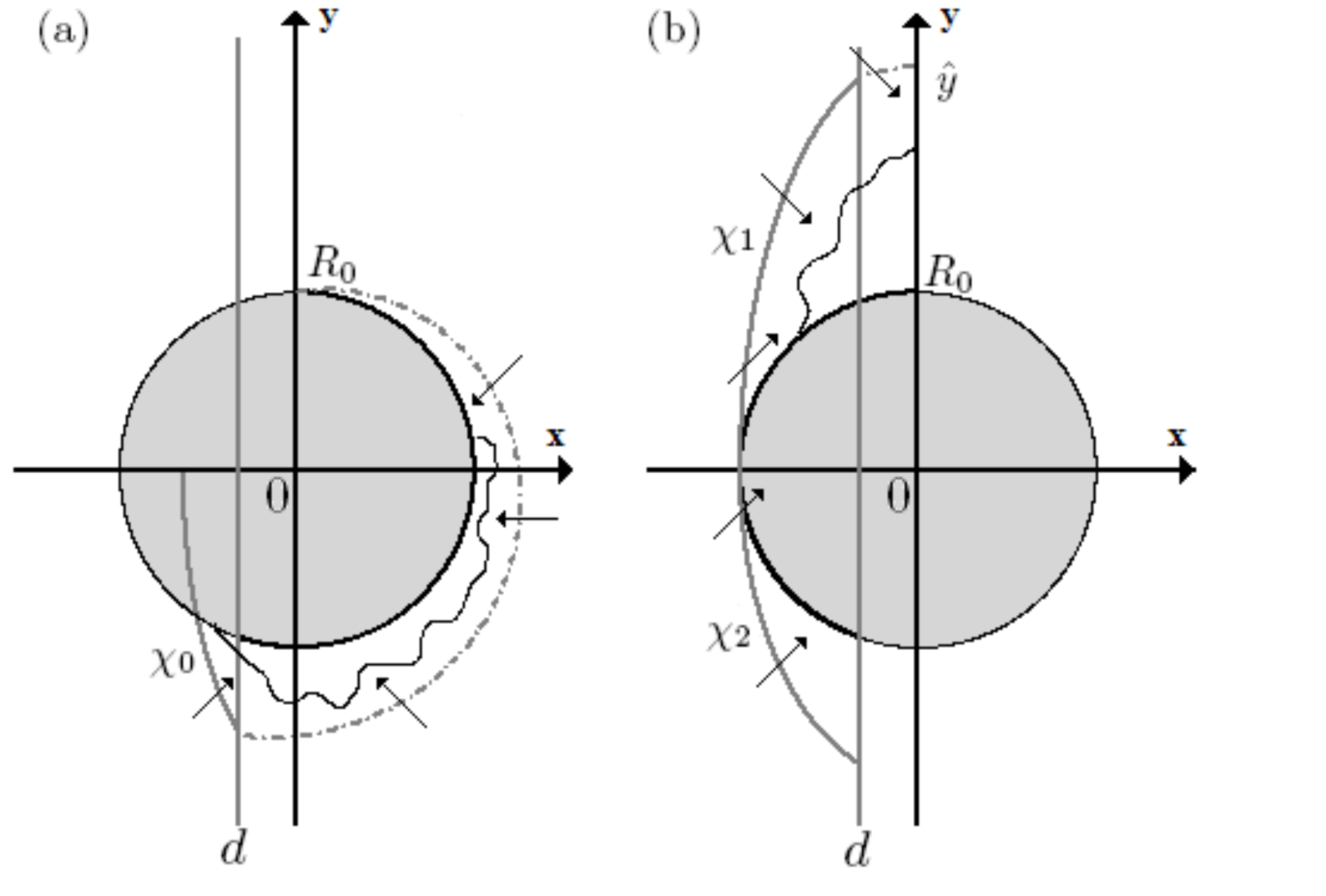, width = 12 cm}}
\caption{(a) If $x(\tau_0)>0$ or $x(\tau_0)>d$ with $y(\tau_0)<0$ then the solution re-enters the ball of radius $R_0$ before exiting the bigger ball of radius $\mathcal R(R_0)$, excluding this situation. (b) In the other cases, the level curves $\chi_1$ and $\chi_2$, respectively of the energy functions $H_1$ and $H_2$, drive the solutions, permitting us to find the desired estimate $y(\tau_2) <\hat y$.}
\label{figelastic}
\end{figure}

Moreover, by~\eqref{as-}, we can suppose $R_0$ sufficiently large to have
\begin{equation}\label{inout}
2F(-r)-r^2 >2F(x)-x^2   \quad \text{for every } r\geq R_0 \text { and } x\in(-r,d)  \,.
\end{equation}
In particular, once fixed $r\geq R_0$, we have
$$
H_i(x,y)<H_i(-r,0)\,, \quad i=1,2\,,
$$
for every $(x,y)$ satisfying $x^2+y^2<r^2$ and $x<d$.
In other words, if $R_0$ is sufficiently large, the level subsets $\chi_i$ of the energy functions $H_i$ passing through the point $(-R_0,0)$ do not enter the open ball of radius $R_0$, see Figure~\ref{figelastic}b. By~\eqref{as-}, we can also suppose that
\begin{equation}\label{en1}
2F_2(-R_0)>R_0^2 e^{2(\kappa\pi+a)}\,.
\end{equation}

We can now state the following lemma.

\begin{lemma}\label{elastic}
There exists $\mathcal R(R_0)>R_0$ such that every $T$-periodic solution of~\eqref{planar} such that $x^2(t_0)+y^2(t_0)>\mathcal R(R_0)$ at a certain time $t_0$ is a $R_0$-large solution.
\end{lemma}

\proof
Set $\mathcal R(R_0)=\mathcal L^{N+2}(\hat y)$ with $\hat y =e^a \sqrt{2F_1(-R_0)+d^2})$.
Argue by contradiction and suppose the existence of a $T$-periodic solution $(x,y)$ of~\eqref{planar}, such that, for some instants $\tau_0$ and $\tau_1$ with $\tau_0<\tau_1\leq \tau_0+T$, it satisfies $x^2(\tau_0)+y^2(\tau_0)=R_0^2$, and $x^2(t)+y^2(t)>R_0^2$ for $t\in(\tau_0,\tau_1)$ and $x^2(\tau_1)+y^2(\tau_1)>\mathcal R(R_0)$.

By Lemma~\ref{rotate+}, set $\tau_2>\tau_0$ the smallest instant such that $x(\tau_2)=0$ and $y(\tau_2)>0$. We prove now that $y(\tau_2) <\hat y$.

First of all, it could not be $x(\tau_0)>0$, or $x(\tau_0)>d$ with $y(\tau_0)<0$: the solution would enter the ball of radius $R_0$ too early (see Figure~\ref{figelastic}a). In fact, using Lemma~\ref{lapontheright}, we can find an instant $\tau_3\in(\tau_0,\tau_2)$, with $x(\tau_3)=d$ and $-e^{\kappa\pi+a} R_0<y(\tau_3)<0$. Then, by~\eqref{en1} and the estimates in~\eqref{en<1} and~\eqref{en<2}, we obtain the contradiction: the solution re-enters the ball guided by the level curve of $H_2(x,y)=H_2(d,-e^{\kappa\pi+a} R_0)$, denoted by $\chi_0$ in Figure~\ref{figelastic}a.

The possibility of having $x(\tau_0)\leq d$ with $y(\tau_0)<0$ is avoided by the {\em guiding} level curve $H_2(x,y)=H_2(-R_0,0)$, denoted by $\chi_2$ in Figure~\ref{figelastic}b. So, it remains the case $x(\tau_0)\leq 0$ with $y(\tau_0)\geq 0$. In this situation the {\em guiding} level curve $H_1(x,y)=H_1(-R_0,0)$, denoted by $\chi_1$ in figure~\ref{figelastic}b, controls the solution when $x<d$ and then by the estimate in Lemma~\ref{lapontheright} we obtain $y(\tau_2) <\hat y$.

Now, in the interval $[\tau_2,\tau_1]$ the solution performs a certain number of complete rotations around the origin, which is less than $N+2$ thanks to Lemma~\ref{Ngiri}. Hence, by~\eqref{elle}, $y(\tau_1) <\mathcal L^{N+2}(y(\tau_2)) < \mathcal L^{N+2}(\hat y) = \mathcal R(R_0)$ thus giving a contradiction.
\cvd

\medbreak

Some easy consequences of the previous lemma are the followings.

\begin{remark}\label{valgono}
If a $T$-periodic solution $x$ of~\eqref{homotopy} satisfies $\|x\|_\infty> \mathcal R(R)$, then Lemmas~\ref{Ngiri},~\ref{lapontheright} and Remark~\ref{r1} hold for $(x,y)=(x,x')$ solution of~\eqref{planar}.
\end{remark}

\begin{remark}\label{tonorm}
Suppose to have a sequence $x_n$ of $T$-periodic solutions of~\eqref{planar} such that $\lim_n \max_{[0,T]} (x_n^2(t)+y_n^2(t)) = +\infty$ then $\lim_n \|x_n\|_\infty = +\infty$.
\end{remark}

The proof of this last remark easily follows by noticing that Lemma~\ref{elastic} holds similarly for every $R>R_0$.

\medbreak

Repeating some of the arguments contained in the proof of Lemma~\ref{elastic}, we can see that $x(t_6)>\mathcal M(x(t_3))$ where
$\mathcal M(\upsilon) = F_2^{-1} \left( \frac12 (e^{\kappa \pi + 2a} \, \upsilon^2-d^2)\right)$. In particular
 
$$
\lim_{r\to+\infty} \frac{\mathcal M(r)}{r} = 0\,.
$$
As an immediate consequence, we have the following lemma.
\begin{lemma}\label{mintozero}
Suppose to have a sequence $x_n$ of $T$-periodic solutions to~\eqref{homotopy}, with $\lim_n \|x_n\|_\infty = +\infty$, then
$$
\lim_{n} \frac{\min x_n(t)}{\| x_n\|_\infty} = 0 \,.
$$
\end{lemma}

\subsection{The a priori bound}\label{priori}

The proof of Theorems~\ref{main} and~\ref{main2} is given essentially by the validity of the following proposition.

\begin{proposition}\label{apb}
There exists $R_{good}$ sufficiently large, such that every $T$-periodic solution of~\eqref{planar} satisfies $x^2(t)+y^2(t)<R_{good}$ for every $t\in[0,T]$.
\end{proposition}

\proof
We argue by contradiction and suppose that there exists a sequence of $T$-periodic solutions $(x_n,y_n)$ of~\eqref{planar}, with $\lambda=\lambda_n$, such that $\max_{[0,T]} x^2_n(t)+y^2_n(t) > n^2$. We have immediately, by Remark~\ref{tonorm}, that $\lim_n \|x_n\|_\infty = +\infty$. Let us denote by $\bar t_n$ the point of maximum of $x_n$, i.e. such that $x_n(\bar t_n)=\|x_n\|_\infty$. We can assume, by Lemma~\ref{elastic}, all these functions to be $R_0$-large. In particular, by Lemma~\ref{Ngiri}, all the solutions must perform $N$ or $N+1$ rotations around the origin.

Consider the sequence of normalized functions
$$
v_n = \frac{x_n}{\|x_n\|_\infty} \,,
$$
which are solutions to
\begin{equation}\label{approx}
v_n'' + \frac{g_{\lambda_n}(t,x_n(t))}{\|x_n\|_\infty}=0\,.
\end{equation}
We have, by Lemma~\ref{mintozero}, $v_n\leq 1=v_n(\bar t_n)$ and $\lim_n \min v_n = 0$. Clearly, $v_n'=y_n/\|x_n\|_\infty$ and, by Remark~\ref{r1}, $\|v_n'\|_\infty \leq e^{\kappa\pi/2}$. For this reason, up to subsequence, $v_n$ converges {\em weakly} in $H^1$ and uniformly to a $T$-periodic non-negative function $v$, with $\|v\|_\infty =1$.
Moreover, we can assume that $\lambda_n \to \bar \lambda$ and that all the solutions $v_n$ draw in the phase-plane the same number of rotations around the origin $K \in\{N,N+1\}$.

\medbreak

We can find some instants $t^n_r$ and $s^n_r$ such that the solutions $v_n$, in the phase-plane $(x,y)$, cross the positive $y$ semi-axis at $t^n_r$ and the negative $y$ semi-axis at $s^n_r$, i.e.
\begin{equation}\label{ts}
t_1^n < s_1^n < t_2^n < s_2^n < \cdots < t_{K}^n < s_{K}^n <  t_{K+1}^n = t_1^n + T \,,
\end{equation}
such that, for every $r\in\{1,\ldots,K\}$,
$$
x_n(t) >0 \text{ for every } t\in(t_r^n, s_r^n)\,, 
$$
$$
x_n(t) <0 \text{ for every } t\in(s_r^n, t_{r+1}^n)\,.
$$
Up to subsequences, we can assume that $t_r^n \to \check\xi_r$ and $s_r^n \to \hat\xi_r$ such that
$$
\check\xi_1 \leq \hat\xi_1 \leq \check\xi_2 \leq \hat\xi_2 \leq \cdots \leq\check\xi_{K} \leq \hat\xi_{K} \leq \check\xi_{K+1} = \check\xi_1 +T\,.
$$
By Lemma~\ref{halflap}, we have $\lim_n t_{r+1}^n - s_r^n =0$, then $\hat\xi _r =\check\xi_{r+1}$. Let us simply denote $\xi_r = \hat\xi_r = \check\xi_{r+1}$. Clearly, $v(\xi_r)=0$.
By the estimate in Lemma~\ref{halflap}, we can easily conclude that necessarily
$\xi_{r+1} - \xi_{r} = T/K$. 

Being $\|v\|_\infty=1$ we are sure that there exists an index $r$ such that $v>0$ in the interval $J_r=(\xi_r,\xi_{r+1})$.
Let us state the following claims, which will be proved in Section~\ref{newsec}, for the reader convenience. We emphasize that the proof of these claims is a crucial part of the proof of Theorems~\ref{main} and~\ref{main2}.

\begin{claim}\label{maxint}
Suppose that $v$ is positive in at least one instant of an interval $J_r=(\xi_r,\xi_{r+1})$, then $v(t)>0$ for every $t\in J_{r}$.
\end{claim}

\begin{claim}\label{alwayspos}
If {\bf(H)} holds, then we have $v>0$ in the interval $J_r=(\xi_r,\xi_{r+1})$, for every index $r$. Moreover, the right and left derivatives at $\xi_r$ exist with $-v'(\xi_r^-)=v'(\xi_r^+)$.
\end{claim}

\medbreak

We will now prove that $v$ solves $v'' + \mu_{K} v =0$ for almost every $t$. By the use of some functions with compact support in $J_r$, we can prove (see~\cite{FGsing} for details) that $v\in H^2_{loc}(J_r)\cap C^1(J_r)$ is a {\em weak} solution of $v'' + p(t) v =0$ in any interval $J_r$, where $p(t)$ is such that $\mu_N\leq p(t) \leq \mu_{N+1}$.

We need to show that $p(t) = \mu_{K}$ for almost every $t\in J_r$.
Consider one of the intervals $J_r$ in which $v$ remains positive (Claim~\ref{maxint} guarantees that $v$ remains positive in the whole interval $J_r$). We will simply denote the extremals of $J_r$ with $\alpha$ and $\beta$, i.e. we set $(\alpha,\beta)=J_r$ for the reader convenience. We have $\beta-\alpha=T/K$.
Introducing modified polar coordinates
$$
\begin{cases}
v(t) = \frac{1}{\sqrt{\mu_{K}}} \, \tilde\rho(t) \cos( \tilde\vartheta(t))\\
v'(t) = \tilde\rho(t) \sin( \tilde\vartheta(t))
\end{cases}
$$
we obtain, integrating $-\tilde\vartheta'$ on $[\alpha,\beta]$, if $K=N$
$$
\pi = \sqrt{\mu_N} \int_\alpha^\beta \frac{p(t)v(t)^2 + v'(t)^2}{\mu_{N}v(t)^2 + v'(t)^2} \,dt \geq \sqrt{\mu_{N}} \, \frac{T}{N} = \pi
$$
and if $K=N+1$
$$
\pi = \sqrt{\mu_{N+1}} \int_\alpha^\beta \frac{p(t)v(t)^2 + v'(t)^2}{\mu_{N+1}v(t)^2 + v'(t)^2} \,dt \leq \sqrt{\mu_{N+1}} \, \frac{T}{N+1} = \pi\,,
$$
thus giving us, in both cases, $p(t)=\mu_{K}$ for almost every $t\in[\alpha,\beta]$. In particular, for every $t\in J_r$, if {\bf (H)} holds
\begin{equation}\label{thisformH1}
v(t) = \sin \big(\sqrt{\mu_{N+1}} (t- \xi_r)\big)\,,
\end{equation}
thanks to Claim~\ref{alwayspos}, while, if it does not hold we have only
\begin{equation}\label{thisform}
v(t) = c_r \sin \big(\sqrt{\mu_{N+1}} (t- \xi_r)\big)\,,
\end{equation}
with $c_r\in[0,1]$ and at least one of them is equal to $1$, being $\|v\|_\infty=1$.
Moreover, we have necessarily $\lambda_n \to \bar\lambda=1$.

\medbreak

We still consider the interval $(\alpha,\beta)=J_r$ for a certain index $r$.
The function $v$ is a solution of the Dirichlet problem:
$$
\begin{cases}
v'' + \mu_{K} \, v =0\\
v(\alpha)=0, \quad v(\beta)=0\,.
\end{cases}
$$

Denote by $\pscal \cdot\cdot$ and $\| \cdot \|_2$, respectively, the scalar product and the norm in $L^2(\alpha,\beta)$.
Call $\phi_K$ the solution of the previous Dirichlet problem with $\| \phi_K\|_2=1$ and introduce the projection of $x_n$ and $v_n$ on the eigenspace generated by $\phi_K$:
$$
x_n^0 = \pscal{x_n}{\phi_K} \phi_K
\quad \aand \quad
v_n^0 = \pscal{v_n}{\phi_K} \phi_K \,.
$$
Being $v=\|v\|_2 \phi_{K}$, we have $v_n^0 \to v$ uniformly in $[\alpha,\beta]$ and $v_n^0\geq 0$, for $n$ sufficiently large.


Multiplying equation~\eqref{homotopy} by $v_n^0$ and integrating in the interval $[\alpha,\beta]$ we obtain
$$
\begin{array}{l}
\ds \int_\alpha^\beta g_{\lambda_n}(t,x_n(t)) v_n^0(t) \,dt 
	= -\int_\alpha^\beta (x_n^0)''(t) v^0_n(t)\,dt\\ [3mm]
\ds\hspace{10mm}=  - \int_\alpha^\beta x_n^0(t) (v_n^0)''(t) \,dt
=   \int_\alpha^\beta \mu_{K} x_n^0(t) v_n^0(t) \,dt \\ [3mm]
 \ds\hspace{20mm}=  \int_\alpha^\beta \mu_{K} x_n(t) v_n^0(t) \,dt \,.
 \end{array}
$$

Defining $r_n(t,x)= g_{\lambda_n}(t,x) - \mu_{K} x$ we have
$$
\int_{\alpha}^{\beta} r_n(t,x_n(t)) v_n^0(t) \,dt = 0
$$

and, applying Fatou's lemma,
$$
\int_{\alpha}^{\beta} \limsup_{n\to\infty} r_n(t,x_n(t)) v_n^0(t) \,dt \geq 0
\geq \int_{\alpha}^{\beta} \liminf_{n\to\infty} r_n(t,x_n(t)) v_n^0(t) \,dt 
\,.
$$
It is easy to see that for every $s_0\in(\alpha,\beta)$ it is possible to find $\bar n(s_0)$ such that $x_n(s_0)>0$ for every $n>\bar n(s_0)$.
So, pointwise, for $n$ large enough $r_n(t,x_n(t))=\lambda_n f(t,x_n(t)) + (1-\lambda_n) \mu x_n(t)-\mu_{K} x_n(t)$.
Hence, being $v_n^0 \to v$ and $\lambda_n \to 1$,
we have, if $K=N$
$$
\int_{\alpha}^{\beta} \liminf_{x\to+\infty} [f(t,x)-\mu_{N} x] v(t) \,dt \leq 0
$$
and if $K=N+1$
$$
\int_{\alpha}^{\beta} \limsup_{x\to+\infty} [f(t,x)-\mu_{N+1} x] v(t) \,dt \geq 0\,.
$$
The previous estimates contradict the hypotheses in~\eqref{LLcond1} if $K=N$ or in~\eqref{LLcond2} if $K=N+1$. Notice that, by Claim~\ref{alwayspos}, if {\bf(H)} holds then this reasoning can be repeated for every interval $J_r$ thus obtaining the contradiction being $v$ as in~\eqref{thisformH1} and not as in~\eqref{thisform}.
\cvd

\subsection{Proof of Claims~\ref{maxint} and~\ref{alwayspos}}\label{newsec} 

In this section we prove Claims~\ref{maxint} and~\ref{alwayspos}. We have preferred to postpone their proof because the arguments we will use are totally independent by the rest of the proof of Proposition~\ref{apb}. This section is inspired by some recent results obtained by the second author in~\cite{S10} for impact systems at resonance (see also~\cite{FS3}).
%
%
%
%
%

The functions $v_n=x_n/\|x_n\|_\infty$ solve equation~\eqref{approx}, which we rewrite in a simpler form
$$
v_n'' + h_n(t,v_n) =0\,,
$$
where, for every $n$,
\begin{equation}\label{hgood}
|h_n(t,v)|\leq d (v + 1) \quad\text{ for every } t\in[0,T] \text{ and } v\geq 0 \,,
\end{equation}
for a suitable constant $d>0$.

\begin{remark}\label{C1conv}
Let $[a,b]\subset J_r$ such that $v$ is positive in $[a,b]$.
The sequence $v_n$ $C^1$-converges to $v$ in $[a,b]$.
\end{remark}

\proof
We have already seen that $(v_n)_n$ is bounded in $C^1$, and by~\eqref{hgood} as an immediate consequence we get $|v_n''(t)|\leq |h_n(t,v_n)|\leq 2d$ for every $t\in[a,b]$. So, being $v_n$ bounded in $C^2$ in such a interval, by the Ascoli-Arzel\` a theorem, we have that $v_n$ $C^1$-converges to $v$ in $[a,b]$.
\cvd

\medbreak

We can now prove the first of the two claims.
%
%

\medbreak

\noindent\textit{Proof of Claim~\ref{maxint}.}
Let $\bar s\in J_{r}$ be the point of maximum of $v$ restricted to the interval $J_r$ with $v(\bar s)=\bar v$. Suppose that $v$ vanishes at $s_0\in J_{r}$, and let $U_0$ be a closed neighborhood of $s_0$ contained in $J_{r}$. We assume without loss of generality that $U_0\subset(\bar s,\xi_{r+1})$, the case $U_0\subset(\xi_{r},\bar s)$ follows similarly.
Notice that $v_n'(t)<-e^{-\kappa\pi/2}\bar v/2<0$ in $U_0$ as a consequence of Lemma~\ref{lapontheright} (cf. Remark~\ref{r1}), so that the previous lemma forces $v$ to be negative on a right neighborhood of $s_0$, thus giving us a contradiction.
\cvd

\medbreak

The following lemma gives us the estimates on the left and right derivatives when $v$ vanishes.

\begin{lemma}\label{der}
Suppose that $v$ is positive in the interval $J_r=(\xi_r,\xi_{r+1})$, then the following limits exist
$$
v'(\xi_r^+)=\lim_{t\to\xi_r^+} v'(t)>0  \quad \text{and} \quad v'(\xi_{r+1}^-)=\lim_{t\to\xi_{r+1}^-} v'(t)<0\,.
$$
\end{lemma}
\proof
We will prove only the existence of the second limit, the other case follows similarly.
In the interval $J_r$ the function $v$ has a positive maximum, thus we can assume that $\max_{J_r} v > M$ and $\max_{J_r} v_n > M$ for a suitable constant $M\in(0,1)$ for large indexes $n$. Using Remark~\ref{r1}, we obtain that $-v_n'(s_{r+1}^n)\in [M/c,cM]$, where $c=e^{\kappa\pi/2}$. So, we can assume up to subsequence that $\lim_n -v_n'(s_{r+1}^n) = \bar y>0$. We now prove that $\lim_{t\to\xi_{r+1}^-} -v'(t)= \bar y$.
Fix $\epsilon>0$ and $s\in(0,\epsilon)$. It is possible to find, for every $n$ sufficiently large, that the following inequalities hold:
$$
\begin{array}{l}
s_{r+1}^n > \xi_{r+1} -s > s_{r+1}^n - 2\epsilon \,,
\\[1mm]
|v_n'(\xi_{r+1}-s) - v'(\xi_{r+1}-s) | <\epsilon\,,
\\[1mm]
|v_n'(s_{r+1}^n)+\bar y|<\epsilon\,.
\end{array}
$$
Moreover, by~\eqref{hgood}, for every $\delta>0$
$$
|v_n'(s_{r+1}^n) - v_n'(s_{r+1}^n-\delta)| < 2d\delta\,,
$$
thus giving us that 
$$
|v'(\xi_{r+1}-s)+\bar y| < (4d+2) \epsilon\,.
$$
The previous inequality holds for every $\epsilon>0$ and $s\in(0,\epsilon)$.
The lemma is thus proved.
\cvd

\medbreak

In what follows we study how the validity of hypothesis {\bf (H)} gives more informations on the function $v$.

\medbreak

\begin{lemma}\label{goodreturn}
Assume {\bf (H)}, then for every index $r$,
$$
\lim_n \frac{-v_n'(s^n_r)}{v_n'(t^n_{r+1})} = \lim_n \frac{-x_n'(s^n_r)}{x_n'(t^n_{r+1})} = 1\,.
$$
\end{lemma}

\proof
Fix $r$ and define the interval $\mathcal I_n=(s^n_r,t^n_{r+1})$, whose length tends to zero for $n$ large.
Using the notation introduced in Figure~\ref{figlapt}, by the estimates in Lemma~\ref{lapontheright}, we can obtain
$$
a^-(y_n(t_4))\leq \sqrt{y_n(t_5)^2+d^2} \leq a^+(y_n(t_4)) \,,
$$
where $a^-(\upsilon)=e^{-a(\upsilon)} \upsilon$ and $a^+(\upsilon)=e^{a(\upsilon)} \upsilon$ with
$a(\upsilon)=\kappa \arcsin(d/\upsilon)$. Moreover, by the same argument which gave us~\eqref{elle}, we have
$$
\mathcal T_{2,1}^{\mathcal I_n}(y_n(t_5)) \leq y_n(t_7) \leq \mathcal T_{1,2}^{\mathcal I_n}(y_n(t_5)) \,,
$$
where
$$
\mathcal T_{i,j}^{\mathcal I_n}(\alpha) = \sqrt{2 F_i^{\mathcal I_n}\left( (F_j^{\mathcal I_n})^{-1} (\alpha^2/2) \right)} \,,
$$
with $F_i^{\mathcal I_n}=F_{i,\zeta_n,\tau_n}$ defined as in {\bf (H)}, being $\mathcal I(\tau_n,\zeta_n)=\mathcal I_n$.
Then, again by Lemma~\ref{lapontheright}, we have
$$
a^-\left(\sqrt{y_n(t_7)^2+d^2}\right)\leq y_n(t_8) \leq a^+\left(\sqrt{y_n(t_7)^2+d^2}\right) \,.
$$
Notice that $\lim_{\upsilon\to \infty} a^\pm(\upsilon)=1$, and by {\bf (H)} we have also 
$$
\lim_n \lim_{\alpha\to \infty} \frac{\mathcal T_{i,j}^{\mathcal I_n}(\alpha)}{\alpha}=1\,.
$$
Hence, the desired estimate follows.
\cvd

\medbreak

The previous estimate is the main ingredient we need to prove the following lemma.

\begin{lemma}\label{pos}
Assume {\bf(H)}. Suppose that $v$ is positive for a certain $t_0\in J_r=(\xi_{r},\xi_{r+1})$, then $\xi_r$ and $\xi_{r+1}$ are isolated zeros. Hence, by Claim~\ref{maxint}, as an immediate consequence $v$ is positive in every interval $J_r$.
\end{lemma}

\proof 
We just prove that $\xi_{r+1}$ is an isolated zero.
By the argument presented in the proof of Lemma~\ref{der}, if the left derivative $v'(\xi_{r+1}^-)=-\eta<0$, then we can assume $-v_n'(s_{r+1}^n)>\eta/2$ for $n$ large enough.
Suppose by contradiction that there exists $\ee_0\in(0,\eta/8d)$, with $d$ as in~\eqref{hgood}, such that $v(\xi_{r+1}+\ee_0)=0$. For every $n$ large enough we have 
$|t_{r+2}^n-\xi_{r+1}|<\ee_0/4$ and by Lemma~\ref{goodreturn} $v_n'(t_{r+1}^n)>\eta/2$.
 Being $|v_n''|\leq 2d$ when $v_n$ is positive we can show that if $s<\eta/4d$ then $v_n(t_{r+2}^n+s)>s\, \eta/4$. By construction $\xi_{r+1}+\ee_0=t_{r+2}^n+s_0$ for a certain $s_0\in(\ee_0/2,\eta/4d)$, so that we obtain $v_n(\xi_{r+1}+\ee_0)=v_n(t_{r+2}^n+s_0)>\eta\ee_0/8$ for every $n$ large enough, thus contradicting $v_n\to v$.
\cvd

\medbreak

We can now prove the remaining claim.

\medbreak

\textit{Proof of Claim~\ref{alwayspos}.} The first part of the statement is given by Lem\-ma~\ref{pos}. The estimate on the derivatives easily follows by Lemmas~\ref{der} and~\ref{goodreturn}, remembering that in the proof of Lemma~\ref{der} we have shown that $\lim_n = v_n'(s^n_r) = v'(\xi_r^-)$ and $\lim_n = v_n'(t^n_{r+1}) = v'(\xi_r^+)$.
\cvd

%
%
%

\medbreak

\section{Nonlinearities with a singularity and\\ radially symmetric systems}\label{secrad}

In this section we provide a result of existence of periodic solutions to scalar differential equations with a singularity in the spirit of Theorems~\ref{main} and~\ref{main2}. In particular we consider the differential equation

\begin{equation}\label{nleq2}
x''+ f(t,x)=0 \,,
\end{equation}
where $f: \R \times (0,+\infty) \to \R$ is a continuous function $T$-periodic in the first variable. The nonlinearity $f$ presents a {\em strong} singularity at $x=0$, in the following sense.
\begin{itemize}
\item[{\bf(${\rm A}^0$)}]{\sl
There exist $\delta>0$ and two continuous functions $f_1, f_2:(0,\delta) \to \R$ such that
$$
f_1(x) < f(t,x) < f_2(x) < 0, \quad \forevery t\in\R \aand x\in(0,\delta) \,,
$$
satisfying
$$
\lim_{x\to0^+} f_i(x) = -\infty \quad \aand \quad \int_0^{\delta}  f_i(\xi)\,d\xi = -\infty\,, \quad i=1,2\,.
$$
}
\end{itemize}

We assume that the nonlinearity $f$ has an asymptotically linear growth at infinity, as follows.

\begin{itemize}
\item[{\bf(${\rm A}^\infty$)}]{\sl
There exist a constant $c>0$ and an integer $N>0$ such that
$$
\mu_N x - c \leq f(t,x) \leq \mu_{N+1} x + c\,,
$$
for every $x>1$ and every $t\in[0,T]$.
}
\end{itemize}

The corresponding of Theorem~\ref{main} can be reformulated for the differential equation~\eqref{nleq2} in this way.

\begin{theorem}\label{mainweak2}
Assume {\bf($A^0$)} and {\bf($A^\infty$)} and the Landesman-Lazer conditions
\begin{equation}\label{LLcond1a2}
\int_0^T \liminf_{x\to+\infty} (f(t,x)-\mu_{N} x ) \phi_N(t+\tau) \, dt > 0 \,,
\end{equation}
\begin{equation}\label{LLcond2a2}
\int_0^T \limsup_{x\to+\infty} (f(t,x)-\mu_{N+1} x ) \phi_{N+1}(t+\tau) \, dt < 0 \,,
\end{equation}
where $\phi_j$ is defined as
$$
\phi_j(t) =
\begin{cases}
\sin (\sqrt{\mu_j} t) & t\in[0, T/j]\\
0 & t\in[T/j, T]
\end{cases}
$$
extended by periodicity.
Then, equation~\eqref{nleq2} has at least one $T$-periodic solution.
\end{theorem}

As in the previous section, we can introduce an additional assumption on the behavior of $f$ near zero, in order to obtain a different version of the previous theorem.

\begin{itemize}
\item[{\bf ($\widetilde{\rm H}$)}] {\sl For every $\tau\in[0,T]$ and for every $\zeta>0$, consider the set $\mathcal I(\tau,\zeta)=[\tau-\zeta,\tau+\zeta]$ and the functions
$$
f_{1,\tau,\zeta}(x) = \min_{t\in \mathcal I(\tau,\zeta)} f(t,x) 
\qquad
f_{2,\tau,\zeta}(x) = \max_{t\in \mathcal I(\tau,\zeta)} f(t,x) 
$$
with their primitives $F_{i,\tau,\zeta}(x)=\int_\delta^x f_{i,\tau,\zeta}(\xi) \,d\xi$.
We assume that
$$
\lim_{\zeta \to 0} \left( \lim_{x\to 0} \frac{F_{2,\tau,\zeta}(x)}{F_{1,\tau,\zeta}(x)} \right) = 1
$$
uniformly in $\tau\in[0,T]$.
}
\end{itemize}

Hence, the corresponding of Theorem~\ref{main2} is the following.

\begin{theorem}\label{mainstrong2}
Assume {\bf($A^0$)}, {\bf($A^\infty$)} and {\bf ($\widetilde{\rm H}$)}, and the Landesman-Lazer conditions~\eqref{LLcond1a2} and~\eqref{LLcond2a2}
where, $\phi_j$ is defined as
$$
\phi_j (t)= \left| \sin (\sqrt{\mu_j} t) \right| \,.
$$
Then, equation~\eqref{nleq2} has at least one $T$-periodic solution.
\end{theorem}

Let us here show some nonlinearities satisfying (or not) hypothesis {\bf ($\widetilde{\rm H}$)}, cf. Example~\ref{ex1}. 

\begin{example}\label{ex2}
Suppose that there exists a function $h:(0,+\infty)\to \R$ satisfying 
$$
\lim_{x\to0^+} h(x) = -\infty\,,
$$
such that
$$
0< \liminf_{x\to 0} \frac {f(t,x)}{h(x)} \leq \limsup_{x\to 0} \frac {f(t,x)}{h(x)} < +\infty \,.
$$
Then, {\bf ($\widetilde{\rm H}$)} holds. As a particular case, suppose that $f$ can be split (for $0<x<1$) as 
$f(t,x)=q(t)h(x)+p(t,x)$ with $q(t)>0$ and $\lim\limits_{x\to 0} \frac{p(t,x)}{f(x)}=0$
uniformly in $t$. In particular we can consider nonlinearities not depending on $t$ when $0<x<1$, or nonlinearities as
$f(t,x) = -(1+\sin^2(t)) x^{-5} - x^{-3}$, or $f(t,x)= -x^{-3} - \sin^2(t) x^{-2}$.

Otherwise, if for example $f(t,x) = -x^{-3} - \sin^2(t) x^{-5}$ when $0<x<1$, then $f$ does not satisfies {\bf ($\widetilde{\rm H}$)}.
\end{example}

The previous theorems can be viewed as the generalization of the result provided by del Pino, Man\'asevich and Montero in \cite{dPMM} to nonlinearities {\em near resonance}. Recently an existence result by the introduction of Lazer-Leach conditions has been proved by Wang in \cite{WWW}, and we recall the result obtained by Fonda and Garrione in \cite{FGsing} where the authors provide a Landesman-Lazer condition {\sl on one side}, roughly speaking, with respect to the smaller eigenvalue.
In particular the previous theorems can be viewed as an answer to~\cite[Remark~2.5]{FGsing}.

\subsection{Proof of Theorems~\ref{mainweak2} and~\ref{mainstrong2}}\label{pf2}

The proof of Theorems~\ref{mainweak2} and~\ref{mainstrong2}, follows step by step the proof of Theorems~\ref{main} and~\ref{main2}, with some wise adjustments. Hence, we will provide only a sketch. We refer to~\cite{FTonNA} for detailed computations in this setting.

Let us underline that, up to a rescaling of the $x$ variable, it is not restrictive to assume $\delta=1$ in {\bf ($A^0$)}.
In \cite{FTonNA}, Fonda and Toader provide an {\em a priori bound} to solutions of equation~\eqref{nleq2} when the nonlinearity satisfies ($A^0$) and the {\em nonresonance} condition
$$
\mu_N < \mu_\downarrow \leq \liminf_{x\to+\infty} \frac{f(t,x)}{x} \leq  \limsup_{x\to+\infty} \frac{f(t,x)}{x}\leq \mu_\uparrow < \mu_{N+1} \,.
$$
As a particular case we find the nonlinearity
$$
h(t,x) =
\begin{cases}
f(t,x) & x<1/2\\
(2x-1) \mu x + (2-2x) f(t,x) & 1/2\leq x\leq 1\\
\mu x & x> 1
\end{cases}
$$
with $\mu=(\mu_N+\mu_{N+1})/2$. Arguing as in Section~\ref{secproof}, we can introduce a family of differential equations
\begin{equation}\label{homotopy2}
x'' + g_\lambda(t,x) = 0\,,
\end{equation}
as in~\eqref{homotopy}, and by standard arguments in degree theory, the proof can be easily obtained when we can find an {\em a priori bound} to the solutions of~\eqref{homotopy2}.
Arguing as in Section~\ref{prel} we consider the corresponding system
\begin{equation}\label{planar2}
\begin{cases}
x'=y\\
-y'= g_\lambda(t,x) \,.
\end{cases}
\end{equation}
which is now defined for $(x,y)\in(0,+\infty)\times \R$. We consider the function
$$
\mathcal N(x,y) = \frac{1}{x^2} + x^2 + y^2\,,
$$
so that, as in~\eqref{Rlarge}, we say that
\begin{equation}\label{Nlarge}
(x,y) \text{ is $\mathcal N_0$-large, if } \mathcal N(x,y)>\mathcal N_0 \text{ for every } t\in[0,T]\,.
\end{equation}

All the results contained in Section~\ref{prel} (wisely adjusted) can be reformulated by the study of the phase portrait when $0<x<1$ and when $x>1$. We list some of them for the reader convenience.

\begin{lemma}\label{fullsingrem}
There exists $\mathcal N_0$ sufficiently large such that every $\mathcal N_0$-large solution of~\eqref{planar2}
rotates clockwise around the point $(1,0)$ performing exactly $N$ or $N+1$ rotations.
\end{lemma}

\begin{lemma}\label{bbb}
For every $\ee>0$ there exists $\mathcal N_\ee$ such that every $\mathcal N_\ee$-large solution $(x,y)$ of~\eqref{planar2}, performing a complete rotation around the point $(1,0)$ in the interval $[t_0,t_2]$,  satisfies
$$
t_1-t_0 \in \left(\frac{T}{N+1}-\ee,\frac TN+\ee\right ) \aand t_2-t_1 <\ee\,,
$$
for a certain $t_1\in(t_0,t_2)$, being $x>1$ in the interval $(t_0,t_1)$ and $0<x<1$ in $(t_1,t_2)$.
\end{lemma}

We refer to \cite{FTonNA} for the detailed computation giving us the previous lemmas. We simply underline that the dynamics when $0<x<1$ (respectively when $x>1$) remember the dynamics of the {\em one-sided superlinear} scalar equation previously studied when $x<0$ (resp. when $x>0$). By the construction of some guiding functions we can prove the following estimates. Notice that the use of guiding functions was adopted also in \cite{FTonNA}, by the use of a {\em general method} presented by Fonda and the author in \cite{FS1}.

\begin{lemma}\label{aaa}
There exists $\mathscr N(\mathcal N_0) > \mathcal N_0$ such that every $T$-periodic solution
of~\eqref{planar2} such that $\mathcal N(x(t_0),y(t_0)) > \mathscr N(\mathcal N_0)$ at a certain time $t_0$ is a $\mathcal N_0$-large
solution.
\end{lemma}

\begin{lemma}\label{aaaa}
Suppose to have a sequence $x_n$ of $T$-periodic solutions to~\eqref{planar2} such that $\lim_n \max_{[0,T]} \mathcal N(x_n(t),y_n(t)) = +\infty$ then $\lim_n \|x_n\|_\infty = +\infty$.
\end{lemma}

\medbreak

All the four preceding lemmas are the main ingredients to obtain the desired a priori bound which is given by the next statement.

\begin{proposition}\label{apb2}
There exists $\mathcal N_{good}$ sufficiently large, such that every $T$-periodic
solution of~\eqref{planar2} satisfies $\mathcal N(x(t),y(t)) < \mathcal N_{good}$ for every $t\in[0,T]$.
\end{proposition}

The proof follows the one of Proposition~\ref{apb}: we assume the existence of a sequence of solutions {\em arbitrarily large} in the sense of~\eqref{Nlarge} and we introduce the {normalized sequence} $v_n=x_n/\|x_n\|_\infty$ converging to a certian non-negative function $v$. We can introduce the instants $t^n_r$ and $s^n_r$ as in~\eqref{ts} requiring now that
$x_n(t) > 1$ for every $t\in(t^n_r,s^n_r)$ and $0<x_n(t) < 1$ for every $t\in(s^n_r, t^n_{r+1})$. Similarly, using Lemma~\ref{bbb}, we can obtain the sequence of instants $\xi_r$ such that $v(\xi_r)=0$, 
being $v_n(t^n_r)=v_n(s^n_r)=1/\|x_n\|_\infty \to 0$ for $n\to \infty$. So, whenever we need to consider an interval when $v$ is positive, we can assume the index $n$ sufficiently large to have $x_n>1$ and argue similarly as in Section~\ref{priori}. The analogues of results in Section~\ref{newsec} follows similarly.

\section{Final remarks}\label{final}

We desire now to show an application of Theorems~\ref{mainweak2} and~\ref{mainstrong2} to radially symmetric systems thanks to a general technique introduced in \cite{FTonNA,FTonProc} by Fonda and Toader. We consider the differential equation
\begin{equation}\label{eqrad}
{\bf x}'' + f(t,|{\bf x}|) \frac{{\bf x}}{|{\bf x} |}=0\,,
\end{equation}
where ${\bf x}\in \R^d$  and ${\text{f}}:\R\times(0,+\infty)\to\R$ is a continuous function, $T$-periodic in the first variable.
By the radial symmetry of the equation, every solution of~\eqref{eqrad} is contained in a plane, so we can pass to polar coordinates and consider solutions to the following system
\begin{equation}\label{sist1}
\begin{cases}
\ds \rho'' - \frac{L^2}{\rho^3}+f(t,\rho)=0  \qquad \rho>0 \\[2mm]
\rho^2 \vartheta' = L \,,
\end{cases}
\end{equation}
where $L\in \R$ is the angular momentum. We are interested in the existence of periodic solutions performing a certain number $\nu$ of revolutions around the origin in the time $kT$ and $T$-periodic in the $\rho$ variable, i.e. such that
\begin{equation}\label{percond}
\begin{array}{l}
\ds \rho(t+T)=\rho(t)\,,\\
\ds \vartheta(t+kT)=\vartheta(t)+2\pi\nu\,. 
\end{array}
\end{equation}

Applying the Fonda-Toader {\em general principle for rotating solutions} (cf. \cite[Theorem~2]{FTonProc}), we obtain as a corollary the following theorem, extending to nonlinearities {\em near resonance} the previous result provided in \cite[Theorem~2]{FTonNA} by the same authors.

\begin{theorem}\label{mainrad}
If the nonlinearity $f$ in~\eqref{eqrad} satisfies the hypotheses of Theorem~\ref{mainweak2} (or Theorem~\ref{mainstrong2}) then for every integer $\nu$, there exists an integer $k_\nu$ such that for every integer $k\geq k_\nu$ equation~\eqref{eqrad} has a $kT$-periodic solution $x_{k,\nu}$ which makes exactly $\nu$ revolutions around the origin in the period~$kT$. In particular the corresponding solution of system~\eqref{sist1} satisfies the periodicity conditions~\eqref{percond}. Moreover there exists a constant $R$, independent by the choice of $\nu$, such that $1/R<|{\bf x}_{k,\nu}(t)|<R$ for every $t$ and if $L_{k,\nu}$ denotes the angular momentum of the solution $x_{k,\nu}$, then $\lim_{k\to\infty} L_{k,\nu} = 0$.
\end{theorem}



\newcommand\mybib[8]{{\bibitem{#1} {\sc #2}, {\em #3}, {#4}~{\bf #5} ({#6}), {#7}--{#8}.}} 
\newcommand\mybibb[4]{{\bibitem{#1} {\sc #2}, {\em #3}, {#4}.}} 


%

\bigbreak

\begin{tabular}{l}
Andrea Sfecci\\
Universit\`a Politecnica delle Marche\\
Dipartimento di Ingegneria Industriale e Scienze Matematiche\\
Via Brecce Bianche 12\\
60131 Ancona\\
Italy\\
e-mail: sfecci@dipmat.univpm.it
\end{tabular}

\medbreak

\noindent Mathematics Subject Classification: 34C25

\medbreak
\noindent Keywords: periodic solutions, resonance, superlinear growth, Landesman-Lazer condition, singularity, radially symmetric systems.

\end{document}